\documentclass{amsart}
\usepackage[dvips]{graphics}
\usepackage{epic,eepic}
\usepackage{eucal}
\usepackage{amsmath,amssymb}

\newtheorem{theorem}{Theorem}[section]

\newtheorem{lemma}[theorem]{Lemma}
\newtheorem{corollary}[theorem]{Corollary}
\newtheorem{example}[theorem]{Example}

\newtheorem{conjecture}[theorem]{Conjecture}

\newcommand{\BB}[1]{\mathbb{#1}}

\newcommand{\FRAK}[1]{\mathfrak{#1}}

\renewcommand{\descriptionlabel}[1]%
   {\hspace\labelsep \upshape\bfseries #1}

\newcommand{\Ann}{\operatorname{Ann}}
\newcommand{\QED}{\hspace*{\fill}\raisebox{.54ex}%
   {$\framebox[.35em]{\/}$}\\ \hspace*{4ex}}
\newcommand{\OL}[1]{\overline{#1}}
\newcommand{\UL}[1]{\underline{#1}}
\newcounter{INPUT}
\def\A{\addtocounter{INPUT}{1}\theINPUT}
\def\B{\theINPUT}
\title
[A formula for ideal lattices of commutative rings]
{A formula for ideal lattices of general commutative rings}
\author{Takashi AOKI,\ Shuzo IZUMI,\ Yasuo OHNO,\ Manabu OZAKI}
\date{\today}
\begin{document}
\maketitle
\noindent
\hrulefill
\begin{center}
A{\small BSTRACT}\\[.3cm]
\begin{minipage}{12.5cm}
{\small
\hspace*{2.3ex}Let $S$ be a set of $n$ ideals of a commutative ring $A$ 
and let $G_{\mbox{\scriptsize even}}$ (respectively 
$G_{\mbox{\scriptsize odd}}$) denote the product of 
all the sums of even (respectively odd) number of ideals of $S$. 
If $n\le 6$ the product of $G_{\mbox{\scriptsize even}}$ 
and the intersection of 
all ideals of $S$ is included in $G_{\mbox{\scriptsize odd}}$. 
In the case $A$ is an Noetherian integral domain, this inclusion 
is replaced by equality if and only if $A$ is a Dedekind domain. 
}
\end{minipage}\end{center}\vspace{.5cm}
\begin{center}
\end{center}
\noindent
\hrulefill\\
\section{Introduction}
We know that the product of the greatest common divisor (GCD) 
and the least common denominator (LCM) of 
two natural numbers $a$ and $b$ is the product $ab$. 
Further, it is known that, given a finite set of natural numbers, 
their GCD is expressed in terms of LCMs of its subsets and that 
their LCM is expressed in terms of GCDs of its subsets. 
(see Wolfram \cite{wolf1}, \cite{wolf2}).
These can be generalized to GCDs and LCMs of ideals of 
a Dedekind domain as follows. 

The set of ideals of a commutative ring form a lattice with 
respect to the order of inclusion. 
The GCD and the LCM of a finite number of ideals 
$\FRAK{a}_1,\ldots,\FRAK{a}_n$ are 
defined to be 
the upper bound 
$$
\frak{a}_1\vee\cdots\vee\mathfrak{a}_n=\mathfrak{a}_1+\cdots+\mathfrak{a}_n
$$ 
and 
the lower bound 
$$
\frak{a}_1\wedge\cdots\wedge\frak{a}_n=\frak{a}_1\cap\cdots\cap\frak{a}_n
$$ 
of them. 
Let $S$ be a set of ideals of a Dedekind domain with $n$ 
elements. If $1\le k\le n$, $G(k)$ (respectively $L(k)$) 
denotes the product of all GCDs (respectively LCMs) of subsets 
of $S$ with $k$ elements. Then we have the following equalities
\vspace{.2cm}.\\
{\small$(*)_n$}\hfill
$\displaystyle G(n)
L(2)L(4)\cdots L(2\left\lfloor n/2 \right\rfloor  )
=
L(1)L(3)\cdots L(2\left\lceil n/2 \right\rceil -1),$
\hfill\vspace{.2cm}\\
{\small$(**)_n$}\hfill
$\displaystyle L(n)
G(2)G(4)\cdots G(2\left\lfloor n/2 \right\rfloor  )
=
G(1)G(3)\cdots G(2\left\lceil  n/2 \right\rceil -1).$
\hfill\vspace{.2cm}\\
Conversely, if $A$ is a Noetherian integral domain, 
even one simple equality\vspace{.2cm}\\
{\small$(*)_2=(**)_2$}\hfill
$(\FRAK{a}\wedge\FRAK{b})(\FRAK{a}\vee\FRAK{b})=\FRAK{a}\cdot\FRAK{b}$
\hfill\vspace{.2cm}\\
for any pair $\FRAK{a}$ and $\FRAK{b}$ 
of ideals of a Noetherian integral 
domain $A$ implies that it is a Dedekind domain. 
These are proved in \S\ref{section-dedekind}. 

Both of these equalities fail in the case of general commutative rings 
as seen in \S\ref{example}. 
We found however, if we replace equalities by inclusions, 
one of them remains valid for $n\le 6$ ideals. 
This is our main result:
\vspace{.2cm}

\noindent
{\small$(\dagger)_n$}\hfill$\displaystyle L(n)
G(2)G(4)\cdots G(2\left\lfloor n/2 \right\rfloor  )
\subset
G(1)G(3)\cdots G(2\left\lceil  n/2 \right\rceil -1).$
\hfill\vspace{.2cm}

The number $n$ of the ideals is restricted at present. 
The case of at most five ideals is proved by ordinary mathematical 
reasoning (see \ref{main}) and the case of six ideals using 
a computer program in \S\ref{six}. 
\section{Some combinatorial formulae}\label{combinatorial}
Let $\tilde{T}$ denote the set of all finite list (sequence) of 
elements of a totally ordered set $T$ permitting repetition 
and $\bar{T}$ the quotient of $\tilde{T}$ identifying 
permutated lists. In other words, 
an element of $\bar{T}$ is a \textit{multiset} 
of elements of $T$. In case of an extensional expression of a multiset, 
we use square parenthesis ``$[$" and ``$]$". 
Hence, if $\alpha_i \neq \alpha_j$ $(i \neq j)$, 
the multiset $[\alpha_1,\ldots,\alpha_n]$ 
can be identified with the set $\{\alpha_1,\ldots,\alpha_n\}$. 
Let us put 
$$
N(n):=\{1,2,\ldots,n\}
.$$
For a multiset $[\alpha_1,\ldots,\alpha_n]$ of elements 
of $T$, we define the multisets 
$$
\OL{S}_k:=\OL{S}_k(\alpha_1,\ldots,\alpha_n)
:=
\bigl[
\alpha_{i_1}\vee\cdots\vee\alpha_{i_k}:\ (i_1,\ldots,i_k)\in{N(n)\choose k}
\bigr],
$$$$
\UL{S}_k:=\UL{S}_k(\alpha_1,\ldots,\alpha_n)
:=
\bigl[
\alpha_{i_1}\wedge\cdots\wedge\alpha_{i_k}:\ (i_1,\ldots,i_k)\in{N(n)\choose k}
\bigr],
$$
where $N(n)\choose k$ stands for the family of all subsets 
with $k$ distinct elements of $N(n)$ 
and $\vee$ and $\wedge$ express the supremum and the infimum 
with respect to the total order. 
We can define the join $\cup$ of elements of $\bar{T}$ in an obvious manner. 
\begin{lemma}\label{max-min}
Let $T$ be a totally ordered set. 
For a multiset $[\alpha_1,\ldots,\alpha_n]$ of elements of $T$, 
the following hold. 
\begin{eqnarray*}
\displaystyle
\OL{S}_n\cup\UL{S}_2\cup\cdots\cup\UL{S}_{2\lfloor n/2\rfloor}
&=&
\UL{S}_1\cup\UL{S}_3\cup\cdots\cup\UL{S}_{2\lceil n/2\rceil-1},
\\
\UL{S}_n\cup\OL{S}_2\cup\cdots\cup\OL{S}_{2\lfloor n/2\rfloor}
&=&
\OL{S}_1\cup\OL{S}_3\cup\cdots\cup\OL{S}_{2\lceil n/2\rceil-1}.
\end{eqnarray*}
\end{lemma}
\textit{Proof.} 
First we assume that $\alpha_i$ are distinct. Then we may 
assume that $\alpha_1<\cdots<\alpha_{n}$. 
The term $\alpha_i$ $(1\le i\le n)$ appears 
$\binom{n-i}{j-1}$ times in $\UL{S}_j$. Then 
$$
\delta_{n,i}
+\binom{n-i}{1}+\binom{n-i}{3}+\cdots+\binom{n-i}{2\lceil (n-i)/2\rceil-1}
$$ 
times on the left side of the first equation
and 
$$
\binom{n-i}{0}+\binom{n-i}{2}+\cdots+\binom{n-i}{2\lfloor (n-i)/2\rfloor}
$$
times on the right. These numbers are easily seen to be equal. 
These prove the first equation. 

Next we prove the case when the multiset 
$[\alpha_1,\ldots,\alpha_n]$ contains repetitions. We may assume that 
$$
\alpha_1=\ldots=\alpha_{i_1}
<\alpha_{i_1+1}=\ldots=\alpha_{i_2}
$$$$
<\alpha_{i_2+1}=\ldots=\alpha_{i_p}
<\alpha_{i_p+1}=\ldots=\alpha_n
.$$
Let 
$$
T'=\{
\beta_1<\ldots<\beta_{i_1}
<\beta_{i_1+1}<\ldots<\beta_{i_2}
$$$$
<\beta_{i_2+1}<\ldots<\beta_{i_p}
<\beta_{i_p+1}<\ldots<\beta_n
\}$$
be another totally ordered set. If we define a map 
$\varphi:\ T'\longrightarrow T$ by $\varphi(\beta_i)=\alpha_i$, 
then $\varphi$ commutes with $\vee$ and $\wedge$:
$$
\varphi(\beta_i\vee\beta_j)=\varphi(\beta_i)\vee\varphi(\beta_j),
\quad
\varphi(\beta_i\wedge\beta_j)=\varphi(\beta_i)\wedge\varphi(\beta_j)
.$$
Hence $\varphi$ commutes with $\OL{S}$ and $\UL{S}$ also. 
Then the formulae for $\beta_i$ follow from those for $\alpha_i$, 
which are proved in the above. 

The second equality follows by taking the dual order. 
\QED
\section{Formulae for Dedekind domains}\label{section-dedekind}
Let $A$ be a commutative ring (not necessarily with unity). 
The set of its ideals form a lattice with respect to the order of 
inclusion. 
We can define the \textit{greatest common divisor} (GCD)
and the \textit{least common multiple} (LCM)
for a finite multiset 
$[\FRAK{a}_1,\ldots,\FRAK{a}_n]$ of ideals of $A$ respectively by
$$
GCD(\FRAK{a}_1,\ldots,\FRAK{a}_n)
:=\frak{a}_1\vee\cdots\vee\frak{a}_n
=\frak{a}_1+\cdots+\frak{a}_n,
$$$$
LCM(\FRAK{a}_1,\ldots,\FRAK{a}_n)
:=\frak{a}_1\wedge\cdots\wedge\frak{a}_n
=\frak{a}_1\cap\cdots\cap\frak{a}_n.
$$
Let us put 
$$
G(k):=G(k;\frak{a}_1,\ldots,\frak{a}_n)
=\prod_{(i_1,\ldots,i_k)\in {N(n)\choose k}}
(\frak{a}_{i_1}+\cdots+\frak{a}_{i_k}),
$$$$
L(k):=L(k;\frak{a}_1,\ldots,\frak{a}_n)
=\prod_{(i_1,\ldots,i_k)\in {N(n)\choose k}}
(\frak{a}_{i_1}\cap\cdots\cap\frak{a}_{i_k}).
$$

Let us consider the following equalities.\vspace{.2cm}\\
{\small$(*)_n$}\hfill$\displaystyle G(n)
{L(2)L(4)\cdots L(2\left\lfloor n/2 \right\rfloor  )}
=
{L(1)L(3)\cdots L(2\left\lceil n/2 \right\rceil -1)},$%
\hfill\vspace{.2cm}\\
{\small$(**)_n$}\hfill$\displaystyle L(n)
{G(2)G(4)\cdots G(2\left\lfloor n/2 \right\rfloor  )}
=
{G(1)G(3)\cdots G(2\left\lceil n/2 \right\rceil -1)}.$%
\hfill\vspace{.2cm}

We see the following by \ref{max-min}. 
\begin{theorem}\label{totally}
Let $A$ be a commutative ring. If the ideal lattice of $A$ 
form a totally ordered set, the equalities 
{\small$(*)_n$} and {\small$(**)_n$} hold. 
\end{theorem}
Let $A$ be an integral domain and $K$ the field of all fractions of $A$. 
(An integral domain is always assumed to be with unity in this paper.)
A subset $\frak{a}\subset K$ is called \textit{fractional ideal} if it is an 
$A$-submodule of $K$ and if there exists $r\in A\setminus\{ 0 \}$ such that 
$r\frak{a}\subset A$. If $\FRAK{a}$ is a fractional ideal, 
$A\frak{a}=\frak{a}$. 
If $\frak{a}$ and $\frak{b}$ are fractional ideals of $K$, we can define 
their products to be the set $\frak{a}\frak{b}$ of all finite 
sums of products of elements of $\frak{a}$ and $\frak{b}$. 

An integral domain $A$ is called a 
\textit{Dedekind domain} if the set 
of all the non-zero fractional ideals form a group with respect to this 
multiplication. Any ideal of a Dedekind domain can be expressed as an 
product of prime ideals uniquely up to order (see Lang \cite{lang}). 
For a real number $k$, 
let $\lceil k \rceil$ (resp. $\lfloor k \rfloor$) denote 
       the smallest integer that is greater or equal to 
(resp. the greatest integer that is smaller or equal to) $k$. 

It is known that {\small$(*)_n$} and {\small$(**)_n$} 
hold for all $n\in\BB{N}$ 
for the ring $\BB{Z}$ of rational integers 
(see the site Wolfman Research \cite{wolf1}, \cite{wolf2}). 
We can say a little more.
\begin{theorem}\label{dedekind}
If $A$ is a Noetherian integral domain, the following 
conditions are equivalent. 
\begin{enumerate}
\item
$A$ is a Dedekind domain
\item
The condition {\small$(*)_n$} and {\small$(**)_n$} hold 
for all multisets $[\FRAK{a}_1,\ldots,\FRAK{a}_n]$ of ideals of $A$.
\item
The condition {\small$(*)_2$} $=$ {\small$(**)_2$} holds
for all pairs $[\FRAK{a}_1,\FRAK{a}_2]$ of 
ideals of $A$.
\end{enumerate}
\end{theorem}
\textit{Proof.} 
{\small (2)}$\Longrightarrow${\small(3)}: Trivial. 

{\small (1)}$\Longrightarrow${\small(2)}:  
Suppose that the power of a prime ideal $\frak{p}$ is just $\alpha_i$ 
in the prime product decomposition of $\frak{a}_i$ $(i=1,\ldots,n)$. 
We use the notations $\OL{S}_k$ and $\UL{S}_k$ at the beginning of 
\S\ref{combinatorial}, considering $\BB{R}$ a totally ordered set 
with respect to the ordinary order. 
Then the power of $\frak{p}$ of the left side of the first equality is 
the sum of the multiset
$$
\UL{S}_n(\alpha_1,\ldots,\alpha_{2m})
\cup\OL{S}_2(\alpha_1,\ldots,\alpha_n)
\cup\OL{S}_4(\alpha_1,\ldots,\alpha_n)
\cup\cdots\cup\OL{S}_{2\lfloor n/2\rfloor}(\alpha_1,\ldots,\alpha_n)
$$
and one on the right side is the sum of
$$
\OL{S}_1(\alpha_1,\ldots,\alpha_n)\cup\OL{S}_3(\alpha_1,\ldots,\alpha_n)
\cup\cdots\cup\OL{S}_{2\left\lceil  n/2 \right\rceil -1}
(\alpha_1,\ldots,\alpha_n)
.$$
These are equal by \ref{max-min}. 
This proves {\small$(*)_n$}. The proof of {\small$(**)_n$} is similar. 

{\small (3)}$\Longrightarrow${\small(1)}:  
Let $\FRAK{a}\subsetneqq A$ be a non-zero ideal of $A$. It is enough to 
prove that it is the product of a finite number of prime ideals 
(see Matsumura \cite{matsumura}, {\bf 11.6}). 
Let $\FRAK{p}$ be an associated prime of $A$-module $A/\FRAK{a}$. 
Note that such an associated prime exists because $A$ is Noetherian. 
Then there exists a non-zero element $\bar{x}:=x+\FRAK{a}\in A/\FRAK{a}$ 
such that $\FRAK{p}$ is the annihilator $\Ann_A(\bar{x})$ of $\bar{x}$. 
Since $xy\in\FRAK{a}$ if and only if $y\in\Ann_A(\bar{x})=\FRAK{p}$ 
for $y\in A$, we have $xA\cap\FRAK{a}=x\FRAK{p}$. 
Hence it follows from the assumption {\small (3)} that 
$$
x\FRAK{p}(xA+\FRAK{a})=(xA\cap\FRAK{a})(xA+\FRAK{a})
=x\FRAK{a}.
$$
Since $A$ is an integral domain, we have 
$\FRAK{p}(xA+\FRAK{a})=\FRAK{a}$. If we put $\FRAK{a}_1:=xA+\FRAK{a}$, 
we have $\FRAK{a}=\FRAK{p}\FRAK{a}_1$. 
By the condition $\bar{x}\neq 0$, we have 
$\FRAK{a}\subsetneqq \FRAK{a}_1$. If $\FRAK{a}_1\subsetneqq A$, 
applying the arguments above to the ideal $\FRAK{a}_1$ instead of 
$\FRAK{a}$, we obtain a prime ideal $\FRAK{p}_2$ and an ideal 
$\FRAK{a}_2$ with $\FRAK{a}_1\subset\FRAK{a}_2$ such that 
$\FRAK{a}_1=\FRAK{p}_1\FRAK{a}_2$, which implies 
$\FRAK{a}=\FRAK{p}_1\FRAK{a}_1\FRAK{a}_2$. Continuing this, 
we obtain sequences of ideals 
$\FRAK{a}\subsetneqq\FRAK{a}_1\subsetneqq\FRAK{a}_2\subsetneqq\cdots$ 
and prime ideals 
$\FRAK{p}=\FRAK{p}_0,\ \FRAK{p}_1,\ \FRAK{p}_2,\ldots$ 
such that 
$\FRAK{a}=\FRAK{p}_0\FRAK{p}_1\cdots\FRAK{p}_i\FRAK{a}_{i+1}$ 
for $0\le i\le n-1$. Since $A$ is Noetherian, there exists a number $n\ge 1$ 
such that $\FRAK{a}_n=A$ and 
$\FRAK{a}=\FRAK{p}_0\FRAK{p}_1\cdots\FRAK{p}_{n-1}$.
\QED
\begin{corollary}\label{converse}
If $A$ is a Noetherian factorial (UFD) domain, the following 
conditions are equivalent. 
\begin{enumerate}
\item
$A$ is a principal ideal domain (PID).
\item
$A$ is a Dedekind domain
\item
The condition {\small$(*)_n$} and {\small$(**)_n$} hold for all 
the multiset of ideals.
\item
The condition {\small$(*)_2$} $=$ {\small$(**)_2$} holds 
for all the multiset of two ideals.
\end{enumerate}
\end{corollary}
\textit{Proof.} 
Implication {\small(1)}$\Longleftrightarrow${\small(2)} is obvious 
from Matsumura \cite{matsumura}, 
{\bf 11.6}, {\bf 20.1}. The conditions {\small(2)}, 
{\small(3)} and {\small(4)} are equivalent by \ref{dedekind}. 
\QED
\section{Examples}\label{example}
Here we show that the equalities stated in the previous section 
for ideals of a Dedekind 
domain fail in a general setting using easy examples. 
\begin{example}
{\upshape
\begin{enumerate}
\item
Let us put 
$\FRAK{a}_1:=\langle x^2y \rangle$ and $\FRAK{a}_2:=\langle xy^2 \rangle$ 
in the polynomial ring $\BB{R}[x,y]$. 
Then 
$$
L(2)G(2)=
(\FRAK{a}_1\cap\FRAK{a}_2)(\FRAK{a}_1+\FRAK{a}_2)
=\langle x^4y^3,x^3y^4 \rangle
$$$$
\subsetneqq \langle x^3y^3 \rangle
=G(1)=L(1).
$$
\item
Let us put 
$\FRAK{a}_1:=\langle x \rangle$, 
$\FRAK{a}_2:=\langle y \rangle$,  
$\FRAK{a}_3:=\langle z \rangle$ 
in the polynomial ring $\BB{R}[x,y,z]$. Then we have 
$$
G(3)L(2)=\langle x,y,z \rangle\langle xyz \rangle^2
\subsetneqq
\langle xyz \rangle^2=L(1)L(3)
$$
\item
Let us put 
$\FRAK{a}_1:=\langle x,y \rangle$, 
$\FRAK{a}_2:=\langle y,z \rangle$,  
$\FRAK{a}_3:=\langle z,x \rangle$ 
in the polynomial ring $\BB{R}[x,y,z]$. Then we have 
$$
G(3)L(2)=\langle x,y,z \rangle\langle x,yz \rangle
         \langle y,zx \rangle\langle z,xy \rangle
$$$$
\supsetneqq
\langle x,y \rangle
\langle y,z \rangle
\langle z,x \rangle
\langle xy,yz,zx \rangle
=L(1)L(3),
$$
because $x^2yz \in  G(3)L(2)\setminus L(1)L(3)$. 

\item
Let us put 
$\FRAK{a}_1:=\langle x^2yz \rangle$, 
$\FRAK{a}_2:=\langle xy^2z \rangle$,  
$\FRAK{a}_3:=\langle xyz^2 \rangle$ 
in the polynomial ring $\BB{R}[x,y,z]$. Then we have 
$$
L(3)G(2)
=
(\FRAK{a}_1\cap\FRAK{a}_2\cap\FRAK{a}_3)
(\FRAK{a}_1+\FRAK{a}_2)
(\FRAK{a}_2+\FRAK{a}_3)
(\FRAK{a}_3+\FRAK{a}_1)
$$$$
=
\langle xyz \rangle^2
\cdot
\langle xyz \rangle^3
\cdot
\langle x,y \rangle
\cdot
\langle y,z \rangle
\cdot
\langle z,x \rangle
$$$$
=
\langle xyz \rangle^5
\cdot
\langle x,y \rangle
\cdot
\langle y,z \rangle
\cdot
\langle z,x \rangle
,$$$$
G(1)G(3)
=
(\FRAK{a}_1\FRAK{a}_2\FRAK{a}_3)
(\FRAK{a}_1+\FRAK{a}_2+\FRAK{a}_3)
$$$$
=
\langle xyz \rangle^4
\cdot\langle xyz \rangle
\cdot\langle x,y,z \rangle
=
\langle xyz \rangle^5
\cdot\langle x,y,z \rangle
\vspace{.2cm}
$$
\noindent
Hence $L(3)G(2)\subsetneqq G(1)G(3)$. 
\end{enumerate}
}
\end{example}
Dedekind domains are Noetherian integral domains. 
Let us see that integrality and Noetherianity are indispensable 
in \ref{dedekind}. 
First we show an example of a ring with zero-divisors which 
satisfy the condition {\small$(*)_n$} and {\small$(**)_n$}. 
\begin{example}
\upshape{
The commutative ring $\BB{Z}_2\times\BB{Z}_2$  
is a Noetherian ring but not an integral domain. There 
are 4 ideals. It is easy to confirm the formulae 
{\small$(*)_n$} and {\small$(**)_n$}. 
}
\end{example}
Next we show a non-Noetherian local integral domain which 
satisfy the condition {\small$(*)_2$} and {\small$(**)_2$}. 
\begin{example}
\upshape{
Let $A:=k[x_1,x_2,\ldots]/\langle x_2^2-x_1,x_3^2-x_2,\ldots \rangle$ 
be the quotient ring of a polynomial ring in a countable number of 
variables over field $k$. This is an integral domain. 
Let $\Bar{f}\in A$ denote the equivalence class 
of $f\in k[x_1,x_2,\ldots]$ of $f$ modulo 
$\langle x_2^2-x_1,x_3^2-x_2,\ldots \rangle$. 
Then $A$ is a local ring with the maximal ideal 
$\FRAK{m}:=\langle \bar{x}_1,\bar{x}_2,\ldots \rangle$. 
Take the localization $B:=A_{\FRAK{m}}$. 
Each non-zero ideal of $B$ is either of the forms
\begin{center}
$I_+(a):=\langle x_n^l:\ n, l\in\BB{N},\ l>2^{n-1}a\rangle$ 
for some $a\in\BB{R}$ with $a\ge 0$
\end{center}
or
\begin{center}
$I(l/2^{n-1}):=\langle x_n^l\rangle$ for some $n\in\BB{N}$ 
and for some $l\in\BB{Z}$ with $l\ge 0$.
\end{center}
The ideal $I_+(a)$ is not finitely generated. 
Hence $B$ is not Noetherian. 
It is obvious that 
$$
a>b \Longrightarrow I(a)\subsetneqq I(b),
$$
and 
$$
I_+(l/2^{n-1})\subsetneqq I(l/2^{n-1})
\quad (n\in\BB{N},\ l\in\BB{Z},\ l\ge 0)
.$$ 
Let $\BB{R}\times \{0,\epsilon\}$ be the totally ordered 
set defined as the direct product, with the lexicographical order, 
of $\BB{R}$ with usual order and $\{0,\epsilon\}$ with $0<\epsilon$. 
Then the set of ideals of $B$ with the order of inclusion 
is isomorphic to a subset of $\BB{R}\times \{0,\epsilon\}$ 
and hence it is totally ordered and the equality 
{\small$(*)_n$} and {\small$(**)_n$} hold by \ref{totally}. 
}
\end{example}
\section{Inclusion formula for general commutative rings}
In this section we assume that $A$ is a commutative ring 
(not necessarily with unity). We have the following. 
\begin{theorem}\label{main}
Suppose that $n\le 6$. 
Take a multiset 
$[\FRAK{a}_1,\ldots,\FRAK{a}_n]$ of ideals of $A$. 
Then we have the following.\\
{\small$(\dagger)_n$}\hfill$\displaystyle
L(n)G(2)G(4)\cdots G(2\left\lfloor{n}/{2}\right\rfloor)
\subset
{G(1)G(3)\cdots G(2\left\lceil {n}/{2}\right\rceil -1)}$
\hfill
\end{theorem}
\textit{Proof.} 
In the expressions below, dots, underlines and overlines 
have no meaning. They are used only for the sake of description. 

(1) Proof of {\small$(\dagger)_1$} and {\small$(\dagger)_2$} are 
very easy.

(2) Proof of {\small$(\dagger)_3$}: $L(3)G(2)\subset G(1)G(3)$; 

We have only to prove:
$$
(\FRAK{a}\cap\FRAK{b}\cap\FRAK{c})
\cdot(\FRAK{a}+\FRAK{b})(\FRAK{a}+\FRAK{c})(\FRAK{b}+\FRAK{c})
\subset
\FRAK{a}\FRAK{b}\FRAK{c}(\FRAK{a}+\FRAK{b}+\FRAK{c})
$$
By the symmetry, we have only to prove that 
$(\FRAK{a}\cap\FRAK{b}\cap\FRAK{c})\FRAK{a}^2\FRAK{b}$ 
and 
$(\FRAK{a}\cap\FRAK{b}\cap\FRAK{c})\FRAK{a}\FRAK{b}\FRAK{c}$
are contained on the right side. This is trivial. 

(3) Proof of {\small$(\dagger)_4$}: $L(4)G(2)G(4)\subset G(1)G(3)$;

We have only to prove:
$$
(\UL{\FRAK{a}\cap\FRAK{b}\cap\FRAK{c}\cap\FRAK{d}})
(\FRAK{a}+\FRAK{b})(\FRAK{a}+\FRAK{c})(\FRAK{a}+\FRAK{d})
(\dot{\FRAK{b}}+\FRAK{c})(\FRAK{b}+\FRAK{d})
(\UL{\FRAK{c}+\FRAK{d}})
\cdot(\dot{\FRAK{a}}+\FRAK{b}+\FRAK{c}+\FRAK{d})
$$$$
\subset
\dot{\FRAK{a}}\dot{\FRAK{b}}\OL{\FRAK{c}\FRAK{d}}
(\FRAK{a}+\FRAK{b}+\FRAK{c})(\FRAK{a}+\FRAK{b}+\FRAK{d})
(\FRAK{a}+\FRAK{c}+\FRAK{d})(\FRAK{b}+\FRAK{c}+\FRAK{d}).
$$
By symmetry, we may replace 
$\FRAK{a}+\FRAK{b}+\FRAK{c}+\FRAK{d}$ by $\FRAK{a}$. 
Since this replacement is symmetric with respect to $\FRAK{b}$ and 
$\FRAK{c}$, we may replace it by $\FRAK{b}$. It is obvious that 
$$
(\FRAK{a}\cap\FRAK{b}\cap\FRAK{c}\cap\FRAK{d})(\FRAK{c}+\FRAK{d})
\subset
\FRAK{c}\FRAK{d}
$$
(underlined parts). 
Thus the inclusion reduces to the obvious
$$
(\FRAK{a}+\FRAK{b})(\FRAK{a}+\FRAK{c})(\FRAK{a}+\FRAK{d})
(\FRAK{b}+\FRAK{d})
$$$$
\subset
(\FRAK{a}+\FRAK{b}+\FRAK{c})(\FRAK{a}+\FRAK{b}+\FRAK{d})
(\FRAK{a}+\FRAK{c}+\FRAK{d})(\FRAK{b}+\FRAK{c}+\FRAK{d}).
$$

(4) Proof of {\small$(\dagger)_5$}: $L(5)G(2)G(4)\subset G(1)G(3)G(5)$; 

We have only to prove:\\[.2cm]
$
\begin{array}{lc}
(\OL{\FRAK{a}\cap\FRAK{b}\cap\FRAK{c}\cap \FRAK{d}\cap \FRAK{e}})&\\
\cdot(\FRAK{a}+\FRAK{b})(\FRAK{a}+\FRAK{c})(\FRAK{a}+\FRAK{d})
(\FRAK{a}+\FRAK{e})\\
\cdot\UL{(\FRAK{b}+\FRAK{c})(\FRAK{b}+\FRAK{d})(\FRAK{b}+\FRAK{e})}
(\dot{\FRAK{c}}+\FRAK{d})(\hat{\FRAK{c}}+\hat{\FRAK{e}})
(\OL{\FRAK{d}+\FRAK{e}})&\\[.2cm]
\cdot(\dot{\FRAK{a}}+\FRAK{b}+\FRAK{c}+\FRAK{d})
(\FRAK{a}+\FRAK{b}+\FRAK{c}+\FRAK{e})
(\FRAK{a}+\FRAK{b}+\FRAK{d}+\FRAK{e})
(\FRAK{a}+\FRAK{c}+\FRAK{d}+\FRAK{e})&\\
\cdot(\dot{\FRAK{b}}+\FRAK{c}+\FRAK{d}+\FRAK{e})&
\end{array}
$\hspace{\fill}\\
$\hspace{\fill}
\begin{array}{cr}
& \subset\dot{\FRAK{a}}\dot{\FRAK{b}}\dot{\FRAK{c}}\OL{\FRAK{d}\FRAK{e}}\\
& \cdot(\FRAK{a}+\FRAK{b}+\FRAK{c})(\FRAK{a}+\FRAK{b}+\FRAK{d})
(\FRAK{a}+\FRAK{b}+\FRAK{e}) 
(\FRAK{a}+\FRAK{c}+\FRAK{d})(\FRAK{a}+\FRAK{c}+\FRAK{e})\\
& \cdot(\FRAK{a}+\FRAK{d}+\FRAK{e})\UL{(\FRAK{b}+\FRAK{c}+\FRAK{d})
(\FRAK{b}+\FRAK{c}+\FRAK{e})
(\FRAK{b}+\FRAK{d}+\FRAK{e})} 
(\hat{\FRAK{c}}+\hat{\FRAK{d}}+\hat{\FRAK{e}})\\[.2cm]
& \cdot(\FRAK{\FRAK{a}}+\FRAK{b}+\FRAK{c}+\FRAK{d}+\FRAK{e}).
\end{array}
$\\[.2cm]
We may replace the factor $(\FRAK{a}+\FRAK{b}+\FRAK{c}+\FRAK{d})$ 
on the left by $\FRAK{a}$ without loss of generality. 
Since this choice is symmetric with respect to 
$(\FRAK{b}$, $\FRAK{c}$, $\FRAK{d}$, $\FRAK{e})$, 
we may replace the factor $(\FRAK{b}+\FRAK{c}+\FRAK{d}+\FRAK{e})$ 
factor by $\FRAK{b}$ without loss of generality. Similarly  
we may replace $(\FRAK{c}+\FRAK{d})$ by $\FRAK{c}$. 
Thus the dotted $\FRAK{a},\ \FRAK{b},\ \FRAK{c}$ 
on the left are included in the dotted 
factors on the right. The product of the overlined 
(resp. underlined, hatted) parts on the 
left is included in the overlined (resp. underlined, hatted) 
factor on the right. 

Thus we have only to prove:\\[.2cm]
$
\begin{array}{lc}
(\FRAK{a}+\FRAK{b})(\FRAK{a}+\FRAK{c})(\FRAK{a}+\FRAK{d})
(\FRAK{a}+\FRAK{e})&\\
\cdot(\FRAK{a}+\FRAK{b}+\FRAK{c}+\FRAK{e})
(\FRAK{a}+\FRAK{b}+\FRAK{d}+\FRAK{e})
(\FRAK{a}+\FRAK{c}+\FRAK{d}+\FRAK{e})&
\end{array}
$\hspace{\fill}\\
\hspace{\fill}
$
\begin{array}{cr}
& \subset(\FRAK{a}+\FRAK{b}+\FRAK{c})(\FRAK{a}+\FRAK{b}+\FRAK{d})
(\FRAK{a}+\FRAK{b}+\FRAK{e})\\
& \cdot(\FRAK{a}+\FRAK{c}+\FRAK{d})(\FRAK{a}+\FRAK{c}+\FRAK{e}) 
(\FRAK{a}+\FRAK{d}+\FRAK{e})\\
& \cdot(\FRAK{a}+\FRAK{b}+\FRAK{c}+\FRAK{d}+\FRAK{e}).
\end{array}
$
\\
If we put 
$$
\FRAK{b}':=\FRAK{a}+\FRAK{b},\ 
\FRAK{c}':=\FRAK{a}+\FRAK{c},\ 
\FRAK{d}':=\FRAK{a}+\FRAK{d},\ 
\FRAK{e}':=\FRAK{a}+\FRAK{e}
$$
this reduced to
$$
\UL{\FRAK{b}'\FRAK{c}'\FRAK{d}'}\FRAK{e}'
(\FRAK{b}'+\FRAK{c}'+\FRAK{e}')
(\FRAK{b}'+\FRAK{d}'+\FRAK{e}')
(\FRAK{c}'+\FRAK{d}'+\FRAK{e}')
$$$$\subset
\UL{(\FRAK{b}'+\FRAK{c}')(\FRAK{b}'+\FRAK{d}')}
(\FRAK{b}'+\FRAK{e}')
\UL{(\FRAK{c}'+\FRAK{d}')}(\FRAK{c}'+\FRAK{e}') 
(\FRAK{d}'+\FRAK{e}')
\cdot
(\FRAK{b}'+\FRAK{c}'+\FRAK{d}'+\FRAK{e}').
\vspace{.2cm}$$
The inclusion 
$$
\FRAK{b}'\FRAK{c}'\FRAK{d}'
\subset
(\FRAK{b}'+\FRAK{c}')(\FRAK{b}'+\FRAK{d}')(\FRAK{c}'+\FRAK{d}')
$$
is obvious. If we put 
$$
\FRAK{b}'':=\FRAK{b}+\FRAK{e},\ 
\FRAK{c}'':=\FRAK{c}+\FRAK{e},\ 
\FRAK{d}'':=\FRAK{d}+\FRAK{e},\ 
$$
we have only to prove that 
$$\FRAK{e'}
(\FRAK{b}''+\FRAK{c}'')
(\FRAK{b}''+\FRAK{d}'')
(\FRAK{c}''+\FRAK{d}'')
\subset
\FRAK{b''}\FRAK{c''}\FRAK{d''}
(\FRAK{b}''+\FRAK{c}''+\FRAK{d}'').
$$
Since $\FRAK{e'}\subset\FRAK{b}''\cap\FRAK{c}''\cap\FRAK{d}''$, 
this reduces to the case $n=4$. 

(5) Proof of {\small$(\dagger)_6$} is done by computer soft 
\textit{Mathematica}. 
The program is shown in the next section. 
\QED
It is natural to conjecture the following.
\begin{conjecture}
\upshape{
The theorem above holds for all natural numbers $n$. 
}
\end{conjecture}
\section{The proof of the inclusion formula for six ideals.}%
\label{six}
Here we check validity of \ref{main} in the 
case of six ideals using the 
computer program on \textit{Mathematica}. 

Observing the proof of the case of $n\le 5$ in the previous 
section, we notice that we have only to treat the polynomials 
in $\BB{Z}[x_1,\ldots,x_k]$ obtained by considering 
the ideals as variables.

Hence we put  
$$
G(k):=G(k;x_1,\ldots,x_6)
=\prod_{(i_1,\ldots,i_k)\in{N(6)\choose k}}%
(x_{i_1}+\cdots+x_{i_k}).
$$
Let $P$ and $Q_0$ denote the set of all the monomials appearing 
in the expansion of $G(2)G(4)G(6)$ and $G(1)G(3)G(5)$ respectively. 
The problem reduces to validity of the following:\\[.2cm]
\hspace*{2em}
\textit{
For any $m \in P$ there exists $i\in N(6)$ 
such that $x_im \in Q_0\vspace{.2cm}$. 
}
\\
If $Q$ denote the set of all the monomials obtained by dividing 
the elements of $Q_0$ by some $x_i$:
$$
Q:=\{ m/x_i:\ m\in Q_0,\ i\in N(6)\}
,$$ 
the assertion reduces to the simple inclusion $P\subset Q$. 
Of course this reduces further to the membership problem of 
multi-exponents. The program runs as follows. 
``In[\ ]:=" and ``Out[\ ]:=" mean input and output respectively and 
$(a,b,c,d,e,f)=(x_1,x_2,x_3,x_4,x_5,x_6)$. 
\\

{\scriptsize
\noindent
In[\A] := \verb/p = {a, b, c, d, e, f};/
\\
In[\A] := \verb/q1 = Apply[Times, p]/
\\
\qquad\qquad Out[\B] := a b c d e f
\\
In[\A] := \verb/q2 = Product[p[[i]]+p[[j]], {i, 1, 5}, {j, i+1, 6}]/
\\
\qquad\qquad Out[\B] := (a+b) (a+c) (b+c) (a+d) (b+d) (c+d) (a+e)
\\
\qquad\qquad\qquad\quad (b+e) (c+e) (d+e) (a+f) (b+f) (c+f) (d+f) (e+f)
\\
In[\A] := \verb/q3 = Product[p[[i]]+p[[j]]+p[[k]], {i, 1, 4}, {j, i+1, 5}, {k, j+1, 6}]/
\\
\qquad\qquad Out[\B] := (a+b+c) (a+b+d) (a+c+d) (b+c+d) (a+b+e) (a+c+e)
\\
\qquad\qquad\qquad\quad (b+c+e) (a+d+e) (b+d+e) (c+d+e) (a+b+f) (a+c+f) (b+c+f)
\\
\qquad\qquad\qquad\quad (a+d+f) (b+d+f) (c+d+f) (a+e+f) (b+e+f) (c+e+f) (d+e+f)
\\
In[\A] := \verb/q4 = Product[Apply[Plus, Delete[p, {{i},{j}}]], {i, 1, 5}, {j, i+1, 6}]/
\\
\qquad\qquad Out[\B] := (a+b+c+d) (a+b+c+e) (a+b+d+e) (a+c+d+e) (b+c+d+e)
\\
\qquad\qquad\qquad\quad (a+b+c+f) (a+b+d+f) (a+c+d+f) (b+c+d+f) (a+b+e+f)
\\
\qquad\qquad\qquad\quad (a+c+e+f) (b+c+e+f) (a+d+e+f) (b+d+e+f) (a+b+e+f)
\\
In[\A] := \verb/q5 = Product[Apply[Plus, Delete[p, {i}]], {i, 1, 6}]/
\\
\qquad\qquad Out[\B] := (a+b+c+d+e) (a+b+c+d+f) (a+b+c+e+f)
\\
\qquad\qquad\qquad\quad (a+b+d+e+f) (a+c+d+e+f) (b+c+d+e+f)
\\
In[\A] := \verb/q6 = Apply[Plus, p]/
\\
\qquad\qquad Out[\B] := a+b+c+d+e+f
\\
In[\A] := \verb/expo[mon_] := {Exponent[mon, a], Exponent[mon, b],/
\\
\qquad\ \verb/Exponent[mon, c], Exponent[mon, d], Exponent[mon, e], Exponent[mon, f]}/
\\
In[\A] := \verb/lhs = ExpandAll[q2 q4 q6];/
\\
In[\A] := \verb/rhs = ExpandAll[q1 q3 q5];/
\\
In[\A] := \verb/leftlist = Table[expo[lhs[[k]], {k, 1, Length[lhs]}];/
\\
In[\A] := \verb/rightlist = Table[expo[rhs[[k]], {k, 1, Length[rhs]}];/
\\
In[\A] := \verb/ee[j_] := ee[j] = IdentityMatrix[6][[j]]/
\\
In[\A] := \verb/rightlist0 = {}/
\\
In[\A] := \verb/Do[rightlist0 = Union[rightlist0, {rightlist[[k]]-ee[j]}],/
\\
\qquad\verb/{k, 1, Length[rhs]}, {j, 1, 6}]/
\\
In[\A] := \verb/temp = {};/
\\
In[\A] := \verb/Do[If[MembreQ[rightlist0, leftlist[[k]]]] == False,/
\\
\qquad\ \verb/AppendTo[temp, leftlist[[k]]], {k, 1, Length[lhs]}]/
\\
In[\A] := \verb/Length[temp]/
\\
\qquad\qquad Out[\B] := $0$
}\\


The final output ``0" means that $P\subset Q$, proving 
the case $n=6$ of the theorem. 
\QED

\noindent
\hrulefill\\[.3cm]
\noindent
\begin{minipage}{15em}
\begin{flushright}
Department of Mathematics\\
Kinki University\\
Kowakae Higashi-Osaka\\
577-8502, Japan
\end{flushright}
\end{minipage}
\hfill
\begin{minipage}{13em}
{\bf e-mails of authors:}\\
{\sf
aoki@math.kindai.ac.jp\\
izumi@math.kindai.ac.jp\\
ohno@math.kindai.ac.jp\\
ozaki@math.kindai.ac.jp}
\end{minipage}\\
\end{document}